\documentclass[11pt, final, reqno]{amsart}
\usepackage[mathscr]{eucal}

\theoremstyle{plain}
\newtheorem{theorem}{Theorem}
\newtheorem{corollary}[theorem]{Corollary}
\newtheorem{proposition}{Proposition}[section]
\newtheorem{lemma}[proposition]{Lemma}

\theoremstyle{remark}

\theoremstyle{definition}
\newtheorem*{example}{Example}

\newcommand{\Aut}{\operatorname{Aut}}
\newcommand{\cB}{\mathscr{B}}

\newcommand{\Cz}{C_0(X, \mathscr{B}(\mathscr{H}))}

\newcommand{\cH}{\mathscr{H}}
\newcommand{\Iso}{\operatorname{Iso}}
\newcommand{\ot}{\otimes}
\newcommand{\refl}{\operatorname{ref}}

\newcommand{\vp}{\varphi}

\begin{document}

\centerline{\large\textbf{REFLEXIVITY OF THE AUTOMORPHISM AND}}

\centerline{\large\textbf{ISOMETRY GROUPS OF THE SUSPENSION OF
$\cB(\cH)$}}

\vskip 2cm
\centerline{LAJOS MOLN\'AR}

\centerline{Institute of Mathematics}

\centerline{Lajos Kossuth University}

\centerline{4010 Debrecen, P.O.Box 12}

\centerline{Hungary}

\centerline{e-mail: \texttt{molnarl@math.klte.hu}}
\vskip .5cm

\centerline{and}

\vskip .5cm
\centerline{M\'AT\'E GY\H{O}RY}

\centerline{Institute of Mathematics}

\centerline{Lajos Kossuth University}

\centerline{4010 Debrecen, P.O.Box 12}

\centerline{Hungary}

\centerline{e-mail: \texttt{gyorym@math.klte.hu}}

\title{}
\thanks{  This research was supported from the following sources:\\
          1) Joint Hungarian-Slovene research project supported by
             OMFB in Hungary and the Ministry of Science and Technology
             in Slovenia, Reg. No. SLO-2/96,\\
          2) Hungarian National Foundation for Scientific Research
             (OTKA), Grant No. T--016846 F--019322,\\
          3) A grant from the Ministry of Education, Hungary, Reg. No.
             FKFP 0304/1997}
\date{November 25, 1997}
\maketitle
\newpage

\vglue 4cm
\centerline{Running head:}
\vskip .3cm

\centerline{The suspension of $\cB(\cH)$}

\vskip 1.5cm
\centerline{Address of the contact person:}
\vskip .3cm

\centerline{LAJOS MOLN\'AR}

\centerline{Institute of Mathematics}

\centerline{Lajos Kossuth University}

\centerline{4010 Debrecen, P.O.Box 12}

\centerline{Hungary}

\centerline{e-mail: \texttt{molnarl@math.klte.hu}}

\newpage
\vglue 5cm
\centerline{\parbox{11cm}
{{\textsc{Abstract.}} The aim of this paper is to show that the
automorphism and isometry groups of the suspension of $\cB(\cH)$,
$\cH$ being a separable infinite dimensional Hilbert space,
are algebraically reflexive. This means that every local automorphism,
respectively local surjective isometry of $C_0(\mathbb R)\otimes
\cB(\cH)$ is an automorphism, respectively a surjective isometry.}}

\newpage
\section{Introduction and Statement of The Results}

The study of reflexive linear subspaces of the algebra $\cB(\cH)$ of all
bounded linear operators on the Hilbert space $\cH$ represents
one of the most active research areas in operator theory (see
\cite{Had} for a beautiful general view of reflexivity of this kind).
In the last decade, similar questions concerning
certain important sets of transfomations acting on Banach algebras
rather than Hilbert spaces have also attracted attention. The
originators of the research in this direction are Kadison and
Larson. In \cite{Kad}, Kadison studied
local derivations from a von Neumann algebra $\mathcal R$ into a
dual $\mathcal R$-bimodule $\mathcal M$. A continuous linear map from
$\mathcal
R$ into $\mathcal M$ is called a local derivation if it agrees with some
derivation at each point (the derivations possibly differring from point
to point) in the algebra.
This investigation was motivated by the study of Hochschild cohomology
of operator algebras. The main result, Theorem A, in \cite{Kad}
states that in the above setting, every local derivation is a
derivation. Independently, Larson and Sourour proved in \cite{LaSo} that
the same conclusion holds true for local derivations of
$\cB(\mathscr X)$, where $\mathscr X$ is a Banach space. Since then,
a considerable amount of work has been done concerning local derivations
of various algebras. See, for example, \cite{Bre, BrSe1, Cri, Shu,
ZhXi}. Besides derivations,
there are at least two other very important classes of transformations
on operator algebras which certainly deserve attention. Namely, the
group of automorphisms and the group of surjective isometries. In
\cite[Some concluding remarks (5), p. 298]{Lar}, Larson initiated
the study of local automorphisms (the
definition should be self-explanatory) of Banach
algebras. In his joint paper with Sourour \cite{LaSo},
it was proved that if $\mathscr X$ is an infinite dimensional Banach
space, then every surjective local automorphism of $\cB(\mathscr X)$ is
an automorphism (see also \cite{BrSe1}).
For a separable infinite dimensional Hilbert
space $\cH$, it was shown in \cite{BrSe2} that the above conclusion
holds true without the assumption on surjectivity, i.e. every local
automorphism of $\cB(\cH)$ is an automorphism.

Let us now define our concept of reflexivity. Let $\mathscr X$ be a
Banach space (in fact, in the cases we are interested in, this is
a $C^*$-algebra) and for any subset $\mathscr E\subset \cB(\mathscr X)$
let
\[
\refl_{al} \mathscr E=\{ T\in \cB(\mathscr X) : Tx\in \mathscr Ex
\text{ for all } x\in \mathscr X\}
\]
and
\[
\refl_{to} \mathscr E=\{ T\in \cB(\mathscr X) : Tx\in \overline{\mathscr
E x} \text{ for all } x\in \mathscr X\},
\]
where bar denotes norm-closure. The collection $\mathscr E$ of
transformations is called algebraically reflexive if $\refl_{al}\mathscr
E=\mathscr E$. Similarly, $\mathscr E$ is said to be topologically
reflexive if $\refl_{to} \mathscr  E=\mathscr E$. In this terminology,
the main result in
\cite{BrSe2} can be reformulated by saying that the automorphism group
of $\cB(\cH)$ is algebraically reflexive. Similarly,
Theorem 1.2 in \cite{LaSo} states that the Lie algebra of all
generalized derivations on $\cB(\mathscr X)$ is algebraically reflexive.
Obviously, the topological reflexivity is a stronger property than the
algebraic reflexivity. Among the previously mentioned papers, there is
only one which concerns topological reflexivity. Namely,
Corollary 2 in \cite{Shu} asserts that the derivation algebra of any
$C^*$-algebra is topologically reflexive. Hence, not only the local
derivations are derivations in this case, but every bounded linear
map which agrees with the limit of some sequence of
derivations at each point, is a derivation.

As for the automorphism groups of $C^*$-algebras, such a general result
as in \cite{Shu} does not hold true. If $\mathcal A$ is a Banach
algebra, then denote by $\Aut(\mathcal A)$ and $\Iso(\mathcal A)$
the group of automorphisms (i.e. multiplicative linear bijections) and
the group of surjective linear isometries of $\mathcal A$, respectively.
Now, if $X$ is an uncountable
discrete topological space, then it is not difficult to verify that
the groups $\Aut(C_0(X))$ and $\Iso(C_0(X))$ of the $C^*$-algebra
$C_0(X)$ of all continuous complex valued functions on $X$ vanishing
at infinity are not algebraically
reflexive. Concerning topological reflexivity, there are even von
Neumann algebras whose automorphism
and isometry groups are not topologically reflexive.
For example, the infinite dimensional commutative von
Neumann algebras acting on a separable Hilbert space have this
nonreflexivity property as it was shown in \cite{BaMo}.
Let us now mention some positive results.
In \cite{MolStud1} we proved that if $\cH$ is a separable infinite
dimensional Hilbert space, then  $\Aut(\cB(\cH))$ and $\Iso(\cB(\cH))$
are topologically reflexive.
In \cite{MolLond} we studied the reflexivity of the automorphism and
isometry groups of $C^*$-algebras in the
famous Brown-Douglas-Fillmore theory, i.e.
the extensions of the $C^*$-algebra of all compact operators
on $\cH$ by commutative separable unital $C^*$-algebras. We proved
there that the groups $\Aut$ and $\Iso$ are
algebraically reflexive in the case of every such extension, but, for
example, in the probably most important case of
extensions by $C(\mathbb T)$ ($\mathbb T$ is the perimeter of the
unit disc), our groups are not topologically reflexive. This result
seems to be surprising even in the case of the Toeplitz extension.

In this present paper we study our reflexivity problem
for the suspension of $\cB(\cH)$. The suspension
$S\mathcal A$ of a $C^*$-algebra $\mathcal A$ is the tensor
product $C_0(\mathbb R)\ot \mathcal A$ which is well-known to be
isomorphic to
the $C^*$-algebra $C_0(\mathbb R, \mathcal A)$ of all continuous
functions from $\mathbb R$ to $\mathcal
A$ which vanish at infinity. The suspension plays very important role
in K-theory since the $K_1$-group of $\mathcal A$ is the $K_0$-group of
$S\mathcal A$. In Corollary~\ref{C:suspref} below we obtain that the
automorphism and isometry groups of the supsension of $\cB(\cH)$ are
algebraically reflexive. In fact, in what follows we consider
more general $C^*$-algebras of the form $C_0(X)\ot
\cB(\cH)\cong C_0(X, \cB(\cH))$, where $X$ is a locally compact
Hausdorff space.

Turning to the results of this paper,
in our fisrt theorem we describe the general form of the elements of
$\Aut( C_0(X, \cB(\cH)))$ and $\Iso( C_0(X,\cB(\cH)))$.
From now on, let $\cH$ stand for an infinite dimensional separable
Hilbert space.

\begin{theorem}\label{T:formautiso}
Let $X$ be a locally compact Hausdorff space. A linear map $\Phi: C_0(X,
\cB(\cH)) \to C_0(X, \cB(\cH))$ is an automorphism if and only if
there exist a function $\tau: X\to \Aut(\cB(\cH))$ and a
bijection $\vp :X \to X$ so that
    \begin{equation}\label{E:formautoc0}
    \Phi(f)(x)=[\tau(x)](f(\vp(x))) \qquad (f\in C_0(X, \cB(\cH)), x\in
    X).
    \end{equation}
Similarly, a linear map $\Phi: C_0(X, \cB(\cH)) \to C_0(X, \cB(\cH))$
is a surjective isometry if and only if
there exist a function $\tau: X\to \Iso(\cB(\cH))$ and a
bijection $\vp :X \to X$ so that $\Phi$ is of the form
\eqref{E:formautoc0}.

Moreover, if the linear map $\Phi:\Cz \to \Cz$ is an automorphism,
respectively a surjective isometry, then for the maps $\tau, \vp$
appearing in \eqref{E:formautoc0} we obtain that $x\mapsto
\tau(x), x\mapsto \tau(x)^{-1}$ are strongly continuous and that
$\vp:X \to X$ is a homeomorphism.
\end{theorem}

In the following two results we show that the
algebraic reflexivity of our groups in the case of $C_0(X)$ implies the
algebraic reflexivity of $\Aut(C_0(X)\ot \cB(\cH))$ and
$\Iso(C_0(X)\ot \cB(\cH))$.

\begin{theorem}\label{T:suspref}
Let $X$ be a locally compact Hausdorff space.
If the automorphism group of $C_0(X)$ is algebraically reflexive, then
so is the automorphism group of $C_0(X, \cB(\cH))$.
\end{theorem}

\begin{theorem}\label{T:convref}
Let $X$ be a $\sigma$-compact locally compact Hausdorff space.
If the isometry group of $C_0(X)$ is algebraically reflexive, then so is
the isometry group of $C_0(X, \cB(\cH))$.
\end{theorem}

To obtain the algebraic reflexivity of the automorphism and isometry
groups of the suspension of $\cB(\cH)$ we prove the following assertion.

\begin{theorem}\label{T:funcref}
Let $\Omega \subset \mathbb R^n$ be an open convex set.
The automorphism and isometry groups of $C_0(\Omega )$ are algebraically
reflexive.
\end{theorem}

The proof of this result will show how difficult it might be to treat
our reflexivity problem for tensor products of general $C^*$-algebras
or even for the suspension of any $C^*$-algebra with algebraically
reflexive automorphism and isometry groups.

Finally, we arrive at the statement announced in the abstract.

\begin{corollary}\label{C:suspref}
The automorphism and isometry groups of the suspension
of $\cB(\cH)$ are algebraically reflexive.
\end{corollary}

As for the natural question of whether the groups above are
topologically reflexive, we have the immediate negative answer as
follows.

\begin{example}
Let $(\vp_n)$ be a sequence of homeomorphisms of $\mathbb R$ which
converges uniformly to a noninjective function $\vp$. Define
linear maps $\Phi_n, \Phi$ on $C_0(\mathbb R, \cB(\cH))$ by
     \[
     \Phi_n(f)=f\circ \vp_n \quad \text{and} \quad
     \Phi(f)=f\circ \vp \qquad (f\in C_0(\mathbb R, \cB(\cH)), n\in
     \mathbb N).
     \]
Then $\Phi_n$ is an isometric automorphism of $C_0(\mathbb R,
\cB(\cH))$, the sequence $(\Phi_n(f))$ converges to $\Phi(f)$ for every
$f\in C_0(\mathbb R, \cB(\cH))$ but $\Phi$ is not surjective.
\end{example}

\section{Proofs}

We begin with the following lemma on a characterization of certain
closed ideals in $C_0(X, \cB(\cH))$.

\begin{lemma}\label{L:ideal}
Let $X$ be a locally compact Hausdorff space.
A closed ideal $\mathscr I$ in $C_0(X, \cB(\cH))$ is of the form
     \[
     \mathscr I=\mathscr I_{x_0}=\{ f\in C_0(X, \cB(\cH))\, :\,
     f(x_0)=0\}
     \]
for some point $x_0\in X$ if and only if $\mathscr I$ is a
proper subset of a maximal ideal $\mathscr{I}_m$ in $C_0(X, \cB(\cH))$,
there is no closed ideal properly inbetween $\mathscr{I}$ and
$\mathscr{I}_m$, and $\mathscr I$ is not
the intersection of two different maximal ideals in $C_0(X, \cB(\cH))$.
\end{lemma}

\begin{proof}
The structure of closed ideals in Banach
algebras of vector valued functions is well-known. See, for
example, \cite[Remark on p. 342]{Nai}. Using this result,
$\mathscr{I}$ is a closed ideal in $C_0(X, \cB(\cH))$ if and only if
it is of the form
     \begin{equation*}
     \mathscr I=\{ f\in C_0(X, \cB(\cH))\, :\, f(x) \in
     \mathcal I_x \},
     \end{equation*}
where every $\mathcal I_x$ is a closed ideal of $\cB(\cH)$, i.e. by the
separability of $\cH$, every $\mathcal I_x$ is either $\{ 0\}$ or
$\mathscr{C}(\cH)$ or $\cB(\cH)$.
By the help of Uryson's lemma on the construction of continuous
functions on $X$ with compact support, one can
readily verify that the maximal ideals in $C_0(X, \cB(\cH))$ are exactly
those ideals which are of the form
     \begin{equation*}\label{E:idform}
     \mathscr I=\{ f\in C_0(X, \cB(\cH))\, :\, f(x_0) \in
     \mathscr{C}(\cH) \}
     \end{equation*}
for some point $x_0\in X$. Now, the statement of the lemma follows
quite easily.
\end{proof}

\begin{proof}[Proof of Theorem~\ref{T:formautiso}]
We begin with the proof of the statement on isometries. Let $\Phi$ be a
surjective linear isometry of $C_0(X, \cB(\cH))$. As a consequence of
a deep result due to Kaup (see, for example, \cite{DFR}) we obtain
that every surjective linear
isometry $\phi$ between $C^*$-algebras $\mathcal A$ and $\mathcal B$
has a certain algebraic property, namely $\phi$ is a triple isomorphism,
i.e. it satisfies the equality
    \[
    \phi(ab^*c)+\phi(cb^*a)=\phi(a)\phi(b)^*\phi(c)+
    \phi(c)\phi(b)^*\phi(a)
    \]
for every $a,b,c\in \mathcal A$. This implies that $\phi$ preserves the
closed ideals in both directions. Indeed, if $\mathcal I\subset
\mathcal A$ is a closed ideal, then we have
    \[
    \phi(a)\phi(b)^*\phi(c)+ \phi(c)\phi(b)^*\phi(a)
    \in \phi(\mathcal I) \qquad (a,c\in \mathcal A, b\in \mathcal I).
    \]
Let $\mathcal I'=\phi(\mathcal I)$. We obtain that
$a{\mathcal I'}^*b+b{\mathcal I'}^*a\in \mathcal I'$ $(a,b \in
\mathcal B)$.
Since $\mathcal I'$ is a closed linear subspace of $\mathcal B$, if
$b$ runs through an approximate identity, we deduce
    \begin{equation}\label{E:ideal1}
    a{\mathcal I'}^*+{\mathcal I'}^*a\in \mathcal I' \qquad (a \in
    \mathcal B).
    \end{equation}
If now $a$ runs through an approximate identity, then we have
    \begin{equation}\label{E:ideal2}
    {\mathcal I'}^*\subset \mathcal I'.
    \end{equation}
We infer from \eqref{E:ideal1} and \eqref{E:ideal2} that
$a\mathcal I'+\mathcal I' a\subset \mathcal I'$ $(a\in \mathcal B)$,
i.e. $\mathcal I'$
is a closed Jordan ideal of $\mathcal B$. It is well-known that in
the case
of $C^*$-algebras, every closed Jordan ideal is an (associative) ideal
(see, for example, \cite[5.3. Theorem]{CiYo}) and hence the same is
true for $\mathcal I'$.

By Lemma~\ref{L:ideal} we infer that our map $\Phi$ preserves the
ideals
     \[
     \mathscr I_{x}=\{ f\in C_0(X, \cB(\cH))\, :\, f(x)=0\}
     \qquad (x\in X)
     \]
in both directions. This gives us that there exists a bijection
$\vp:X\to X$ for which
     \begin{equation}\label{E:welldef}
     \Phi(f)(x)=0 \Longleftrightarrow f(\vp(x))=0
     \end{equation}
holds true for every $f\in C_0(X, \cB(\cH))$ and $x\in X$.
For any $x\in X$, let us define $\tau(x)$ by the formula
     \begin{equation}\label{E:bibi}
     [\tau(x)](f(\vp(x)))=\Phi(f)(x) \qquad (f\in C_0(X, \cB(\cH))).
     \end{equation}
Because of \eqref{E:welldef} we obtain that $\tau(x)$ is a well-defined
injective linear map on $\cB(\cH)$. Since $\Phi$ is surjective, we have
the surjectivity of $\tau(x)$. Now, we compute
     \begin{gather*}
     [\tau(x)](f(\vp(x))g(\vp(x))^*f(\vp(x)))=\Phi(fg^*f)(x)=\\
     \Phi(f)(x)\Phi(g)(x)^*\Phi(f)(x)=
     [\tau(x)](f(\vp(x)))([\tau(x)](g(\vp(x))))^*[\tau(x)](f(\vp(x)))
     \end{gather*}
for every $f,g\in C_0(X, \cB(\cH))$. This implies that $\tau(x)$ is a
triple automorphism  of $\cB(\cH)$. Since the triple homomorphisms
preserve the partial isometries and every operator
with norm less than 1 is the average of unitaries, it follows that
$\tau(x)$ is a contraction. Applying the same argument to the inverse
of $\tau(x)$, we obtain that $\tau(x)\in \Iso(\cB(\cH))$.
This proves that $\Phi$ is of the form \eqref{E:formautoc0} given in
the statement of our theorem.

Let now $\Phi:C_0(X, \cB(\cH))\to C_0(X, \cB(\cH))$ be a linear map
of the form
    \begin{equation}\label{E:bugybugy}
    \Phi(f)(x)=[\tau(x)](f(\vp(x))) \qquad (f\in C_0(X, \cB(\cH)), x\in
    X),
    \end{equation}
where $\tau:X \to \Iso(\cB(\cH))$ and $\vp:X\to X$ is a bijection.
The function $\vp$ is continuous. Indeed, this follows easily
from the equality $\|f(\vp(x))\|=\|\Phi(f)(x)\|$
and from Uryson's lemma.
To see the strong continuity of $\tau:X \to \Iso(\cB(\cH))$, let
$(x_\alpha)$ be a net in $X$ converging to $x\in X$. Let
$y_\alpha=\vp(x_\alpha), y=\vp(x)$. We may suppose that every
$y_\alpha$ belongs to a fixed compact neighbourhood of $y$. If $f\in
C_0(X)$ is identically 1 on this neighbourhood, then for every operator
$A\in \cB(\cH)$ we have
      \begin{gather*}
      [\tau(x_\alpha)](A)=[\tau(x_\alpha)](f(\vp(x_\alpha))A)=
      \Phi(fA)(x_\alpha)\longrightarrow \Phi(fA)(x)=\\
      [\tau(x)](f(\vp(x))A)=[\tau(x)](A).
      \end{gather*}
Next, from the equality
      \begin{gather*}
      \| [\tau(x_\alpha)^{-1}](A)-[\tau(x)^{-1}](A)\|=
      \| [\tau(x_\alpha)^{-1}\tau(x)\tau(x)^{-1}](A)-
         [\tau(x)^{-1}](A)\|=\\
      \| [\tau(x)]([\tau(x)^{-1}](A))-
         [\tau(x_\alpha)]([\tau(x)^{-1}](A))\|
      \end{gather*}
we get the strong continuity of the map $x\mapsto \tau(x)^{-1}$.
We prove that $\vp^{-1}$ is also continuous. Since $\Phi$
maps into $C_0(X, \cB(\cH))$, it is quite easy to see from
\eqref{E:bugybugy} that $f\circ \vp\in C_0(X)$ holds true for every
$f\in C_0(X)$. If $K\subset X$ is an arbitrary compact set and $f\in
C_0(X)$ is a function which is identically 1 on $K$, then it follows
from $f\circ \vp\in C_0(X)$
that there exists a compact set $K'\subset X$ for which $\vp(x)\in
K^c$ holds true for all $x\in {K'}^c$. Thus, we have $K\subset \vp(K')$.
Let $(x_\alpha)$ be a net in $X$ such that $(\vp(x_\alpha))$ converges
to some $\vp(x)$. Obviously, we may suppose that every $\vp(x_\alpha)$
belongs to a compact neighbourhood $K$ of $\vp(x)$. By what we have just
seen, there exists
a compact set $K'\subset X$ which contains the net $(x_\alpha)$ and the
point $x$ as well. Since $K'$ is compact, the net
$(x_\alpha)$ has a convergent subnet. Because of the continuity of the
bijection $\vp$, it is easy to see that the limit of this subnet is $x$.
The continuity of $\vp^{-1}$ is now apparent. Finally, one can
verify quite readily that $\Phi$ is a surjective linear isometry of
$C_0(X,\cB(\cH))$.

Let us turn to the proof of our statement concerning automorphisms.
So, let $\Phi$ be an automorphism of $C_0(X,\cB(\cH))$.
Since, as it is well-known, every automorphism of a $C^*$-algebra
is continuous (in fact, its norm equals the norm of its inverse),
one can get the form \eqref{E:formautoc0} in a way very similar to
that was followed in the case of isometries.
Let now $\Phi:C_0(X,\cB(\cH))\to C_0(X,\cB(\cH))$ be a linear map
of the form
    \begin{equation}\label{E:sajo}
    \Phi(f)(x)=[\tau(x)](f(\vp(x))) \qquad (f\in C_0(X, \cB(\cH)), x\in
    X),
    \end{equation}
where $\tau:X \to \Aut(\cB(\cH))$ and $\vp:X\to X$ is a bijection.
We show that $\vp$ is continuous.
Let $(x_\alpha)$ be a net in $X$ converging to $x\in X$.
By \eqref{E:sajo} we have
    \[
    f(\vp(x_\alpha))I=
    [\tau(x_\alpha)](f(\vp(x_\alpha))I)\longrightarrow
    [\tau(x)](f(\vp(x))I)=f(\vp(x))I
    \]
for every $f\in C_0(X)$. Referring to Uryson's lemma again, we infer
that $\vp(x_\alpha)\to \vp(x)$. This verifies the continuity of $\vp$.
We claim that the function $\tau$ is bounded. In fact,
by the principle of uniform boundedness, in the opposite case
we would obtain that
there exists an operator $A\in \cB(\cH)$ for which $[\tau(.)](A)$ is not
bounded. Then there is a sequence $(x_n)$ in $X$ with the property that
$\| [\tau(x_n)](A)\| >n^3$ $(n\in \mathbb N)$. Using Uryson's lemma, it
is an easy task to construct a nonnegative function $f\in C_0(X)$ for
which $f(\vp(x_n))\geq 1/n^2$. Indeed, for every $n\in \mathbb N$ let
$f_n:X
\to [0,1]$ be a continuous function with compact support such that
$f_n(\vp(x_n))=1$ and define $f=\sum_n (1/n^2)f_n$. We have
$\|\Phi(fA)(x_n)\|=\|f(\vp(x_n))[\tau(x_n)](A)\|>n$ $(n\in \mathbb N)$
which contradicts the boundedness of the function $\Phi(fA)$.
The strong continuity of $\tau$ can be proved as it was done in the
case of isometries. Using the inequality
      \begin{gather*}
      \| [\tau(x_\alpha)^{-1}](A)-[\tau(x)^{-1}](A)\|=
      \| [\tau(x_\alpha)^{-1}\tau(x)\tau(x)^{-1}](A)-
         [\tau(x)^{-1}](A)\|\leq \\
      \| \tau(x_\alpha)^{-1}\|
      \| [\tau(x)]([\tau(x)^{-1}](A))-
         [\tau(x_\alpha)]([\tau(x)^{-1}](A))\| =\\
      \| \tau(x_\alpha)\|
      \| [\tau(x)]([\tau(x)^{-1}](A))-
         [\tau(x_\alpha)]([\tau(x)^{-1}](A))\|
      \end{gather*}
and the boundedness of $\tau$, we get the strong continuity of the map
$x\mapsto \tau(x)^{-1}$. The proof can be completed as in
the case of isometries.
\end{proof}

The following two lemmas are needed in the proof of
Theorem~\ref{T:suspref}.

\begin{lemma}\label{L:autcomb}
Let $\tau, \tau_1, \tau_2$ be automorphisms of $\cB(\cH)$ and
let $\lambda \in \mathbb C$, $0\neq \lambda_1, \lambda_2 \in \mathbb C$
be scalars so that
    \[
    \lambda \tau(A)=\lambda_1 \tau_1(A)+\lambda_2 \tau_2(A) \qquad (A\in
    \cB(\cH)).
    \]
Then we have $\tau_1=\tau_2$.
\end{lemma}

\begin{proof}
Since the automorphisms of $\cB(\cH)$ are all spatial (see, for
example, \cite[3.2. Corollary]{Che}),
hence there exist invertible operators $T,T_1,T_2\in \cB(\cH)$ such that
    \begin{equation}\label{E:fugg}
    \lambda TAT^{-1}=\lambda_1 T_1A{T_1}^{-1}+ \lambda_2 T_2A{T_2}^{-1}
    \qquad (A\in \cB(\cH)).
    \end{equation}
It is apparent that if $a,b, x,y, u,v \in X$ and
    \[
    a\ot b=x\ot y+u\ot v,
    \]
then either $\{ x,u\}$ or $\{ y,v\}$ is linearly dependent.
Using this elementary observation and putting $A=x\ot y$ into
\eqref{E:fugg}, we infer that either
$\{ {T_1}x, {T_2}x\}$ is linearly dependent for all $x\in
\cH$ or $\{ {{T_1}^{-1}}^*y, {{T_2}^{-1}}^*y\}$ is linearly dependent
for all $y\in \cH$. In both cases we have the linear
dependence of $\{ T_1, T_2\}$ which results in $\tau_1=\tau_2$.
\end{proof}

In the proof of the next lemma we need the concept of Jordan
homomorphisms. A linear map $\phi$ between algebras $\mathcal A$ and
$\mathcal B$ is called a Jordan homomorphism if it satisfies
     \[
     \phi(A)^2=\phi(A^2) \qquad (A\in \mathcal A).
     \]
If, in addition, $\mathcal A$ and $\mathcal B$ have involutions and
     \[
     \phi(A)^*=\phi(A^*) \qquad (A\in \mathcal A),
     \]
then we say that $\phi$ is a Jordan *-homomorphism.

\begin{lemma}\label{L:locaut}
Let $\Phi:\cB(\cH) \to \cB(\cH)$ be a bounded linear map with the
property that for every $A\in \cB(\cH)$ there exist a
number $\lambda_A\in \mathbb C$ and an automorphism $\tau_A\in
\Aut(\cB(\cH))$
so that $\Phi(A)=\lambda_A \tau_A(A)$. Then there exist a number
$\lambda\in \mathbb C$ and an automorphism $\tau\in \Aut(\cB(\cH))$ such
that $\Phi(A)=\lambda \tau(A)$ $(A\in \cB(\cH))$.
\end{lemma}

\begin{proof}
First suppose that $\Phi(I)=0$. Assume that there exists
a projection $0,I\neq P\in \cB(\cH)$ for which $\Phi(P)\neq 0$. Applying
an appropriate transformation, we may suppose that $\Phi(P)=P$.
Then we have $\Phi(I-P)=-P$. If $\epsilon, \delta$ are
different nonzero numbers, then
by our assumption we infer that $\Phi(\epsilon P+\delta (I-P))$ is a
scalar multiple of an invertible operator which, on the other hand,
equals $(\epsilon-\delta)P$. This clearly implies that
$\epsilon=\delta$, which is
a contradiction. Hence, we obtain that $\Phi(P)=0$ holds true for
every projection $P\in \cB(\cH)$. Using the spectral theorem and the
continuity of $\Phi$, we conclude that $\Phi=0$.

Next suppose that $\Phi(I)\neq 0$. Apparently,
we may assume that $\Phi(I)=I$. By the linearity of $\Phi$,
for an arbitrary projection $0,I\neq P\in \cB(\cH)$ we obtain
    \begin{equation*}\label{E:nondep1}
    I=\Phi(I)=\Phi(P)+\Phi(I-P)=\lambda_{P} Q+\lambda_{I-P} R,
    \end{equation*}
where $Q,R$ are idempotents different from $0,I$. Taking squares on both
sides in the equality
    \[
    I-\lambda_{I-P}R=\lambda_P Q,
    \]
we have
    \[
    I+\lambda_{I-P}^2R-2\lambda_{I-P}R=\lambda_P^2Q.
    \]
But we also have
    \[
    \lambda_{P}(I-\lambda_{I-P}R)=\lambda_P^2 Q.
    \]
Comparing these equalities and using $R\neq 0,I$, we deduce that
$\lambda_P=1$. This means that $\Phi(P)$ is an idempotent. Therefore,
$\Phi$ sends projections to idempotents.
Now, a standard argument shows that $\Phi$ is a
Jordan endomorphism of $\cB(\cH)$ (see, for example, the proof of
\cite[Theorem 2]{MolStud1}). Clearly,
the range of $\Phi$ contains a rank-one operator (e.g. a rank-one
idempotent) and an operator with dense range (e.g. the identity). Using
our former result \cite[Theorem 1]{MolStud1}, we infer that $\Phi$ is
either an automorphism or an antiautomorphism.
This latter concept means that $\Phi$ is a bijective linear map with the
property that $\Phi(AB)=\Phi(B)\Phi(A)$ $(A,B\in \cB(\cH))$. But
$\Phi$ cannot be an antiautomorphism. In fact, in this case we would
obtain that the image $\Phi(S)$ of a unilateral shift $S$ has a right
inverse. But, on the other hand, since $\Phi$ is locally a scalar
multiple of an automorphism of $\cB(\cH)$, it
follows that $\Phi(S)$ is not right invertible. This
contradiction justifies our assertion.
\end{proof}

Before proving Theorem~\ref{T:suspref} we recall that
the automorphisms of the function algebra $C_0(X)$ are of the form
$f\mapsto f\circ \vp$, where $\vp:X\to X$ is a homeomorphism.

\begin{proof}[Proof of Theorem~\ref{T:suspref}]
Let $\Phi:\Cz \to \Cz$ be a local automorphism of $\Cz$,
i.e. $\Phi$ is a bounded linear map which agrees with some automorphism
at each point in $\Cz$.
By Theorem~\ref{T:formautiso}, for every $f\in \Cz$ there exist
a homeomorphism $\vp_f:X\to X$ and a function
$\tau_f: X\to \Aut(\cB(\cH))$ such that
    \begin{equation*}
    \Phi(f)(x)=[\tau_f(x)](f(\vp_f(x))) \qquad (x\in X).
    \end{equation*}
It follows that for every $f\in C_0(X)$ there exists a homeomorphism
$\psi_f:X \to X$ for which $\Phi(fI)=(f\circ \psi_f)I$. Since, by
assumption, the
automorphism group of $C_0(X)$ is reflexive, we obtain that there is a
homeomorphism $\vp:X\to X$ for which
    \begin{equation}\label{E:formcent}
    \Phi(fI)=(f\circ \vp)I \qquad (f\in C_0(X)).
    \end{equation}
Let $f\in C_0(X)$ and $x\in X$. Consider the
linear map $\Psi: A\mapsto \Phi(fA)(x)$ on $\cB(\cH)$. From
the form \eqref{E:formautoc0} of the automorphisms of
$\Cz$ it follows
that $\Psi$ has the property that for every $A\in \cB(\cH)$ there
exist a number $\lambda_A$ and an automorphism $\tau_A\in
\Aut(\cB(\cH))$ such that
    \begin{equation*}\label{E:locform2}
    \Psi(A)=\lambda_A \tau_A(A).
    \end{equation*}
Now, Lemma~\ref{L:locaut} tells us that there exist
functions $\tau_f: X\to \Aut(\cB(\cH))$ and $\lambda_f :X\to \mathbb C$
such that
    \begin{equation*}
    \Phi(fA)(x)=[\tau_f(x)](\lambda_f(x) A) \qquad (f\in C_0(X), A\in
    \cB(\cH), x\in X).
    \end{equation*}
From \eqref{E:formcent} we obtain that $\lambda_f=f\circ \vp$ and
hence we have
    \begin{equation}\label{E:formdepf}
    \Phi(fA)(x)=[\tau_f(x)](f(\vp(x)) A) \qquad (f\in C_0(X), A\in
    \cB(\cH), x\in X).
    \end{equation}
Let $x\in X$ be fixed for a moment. Pick functions $f,g\in C_0(X)$ with
the property that $f(\vp(x)),g(\vp(x))\neq 0$.
Because of linearity we get
    \begin{multline*}
    [\tau_f(x)](f(\vp(x))A)+[\tau_g(x)](g(\vp(x))A)=\\
    \Phi(fA)(x)+\Phi(gA)(x)=\Phi((f+g)A)(x)=\\
    [\tau_{f+g}(x)](f(\vp(x))A+g(\vp(x))A) \quad
    (A\in \cB(\cH)).
    \end{multline*}
Using Lemma~\ref{L:autcomb} we infer that $\tau_f(x)=\tau_g(x)$.
By the formula \eqref{E:formdepf} it follows readily that
there is a function $\tau:X \to \Aut(\cB(\cH))$ for which
    \begin{equation}\label{E:formal}
    \Phi(fA)(x)=[\tau(x)](f(\vp(x)) A)
    \qquad (f\in C_0(X), A\in \cB(\cH), x\in X).
    \end{equation}
Since the linear span of the set of functions $fA$ $(f\in
C_0(X), A\in \cB(\cH))$
is dense in $C_0(X, \cB(\cH))$ (see, for example, \cite[6.4.16.
Lemma]{Mur}), the equality in \eqref{E:formal} gives us that
    \begin{equation*}\label{E:form}
    \Phi(f)(x)=[\tau(x)](f(\vp(x))) \qquad (x\in X)
    \end{equation*}
holds true for every $f\in C_0(X, \cB(\cH))$. By
Theorem~\ref{T:formautiso}, the proof is complete.
\end{proof}

The next lemma that we shall make use in the proof of
Theorem~\ref{T:convref} states that every bounded linear map on
$\cB(\cH)$ which is locally a scalar multiple of a surjective
isometry, equals globally a scalar multiple of a surjective isometry.
For the proof we recall the folk result (in fact this is a
consequence of a theorem of Kadison) that
every surjective linear isometry of $\cB(\cH)$ is either of the form
     \[
     A \longmapsto UAV
     \]
or of the form
     \[
     A \longmapsto UA^{tr}V,
     \]
where $U,V$ are unitary operators and ${}^{tr}$ denotes the transpose
with respect to an arbitrary but fixed complete orthonormal system in
$\cH$. In what follows $\mathscr{P}(\cH)$ and $\mathscr U(\cH)$ denote
the set of all projections and all unitaries on $\cH$, respectively.

\begin{lemma}\label{L:locisom}
Let $\Phi:\cB(\cH) \to \cB(\cH)$ be a bounded linear map with the
property that for every $A\in \cB(\cH)$ there exist a number $\lambda_A
\in \mathbb
C$ and a surjective linear isometry $\tau_A\in \Iso(\cB(\cH))$ so that
$\Phi(A)=\lambda_A \tau_A(A)$. Then there exist a number $\lambda \in
\mathbb C$ and a surjective linear isometry $\tau\in \Iso(\cB(\cH))$ for
which $\Phi(A)=\lambda \tau(A)$ $(A\in \cB(\cH))$.
\end{lemma}

\begin{proof}
Just as in the proof of Lemma~\ref{L:locaut}, first suppose that
$\Phi(I)=0$. Assume that there exists a projection $0,I\neq P\in
\cB(\cH)$ for which $\Phi(P)\neq 0$. Apparently, we
may suppose that $\Phi(P)=P$. Then we have $\Phi(I-P)=-P$. Since for
any different nonzero numbers $\epsilon,\delta \in \mathbb C$, the
operator $\epsilon P+\delta(I-P)$ is invertible, we obtain that
$(\epsilon-\delta)P= \Phi(\epsilon P+\delta (I-P))$ is a scalar
multiple of an invertible operator. But this is a contradiction and
hence we have $\Phi(P)=0$ for every projection $P$. This gives us that
$\Phi=0$.

So, let us suppose that $\Phi(I)\neq 0$. Clearly, we may assume that
$\Phi(I)=I$ and that the constants $\lambda_A$ are all nonnegative.
Let $P\neq 0,I$ be a projection. Let $\lambda,\mu$ be nonnegative
numbers and let $U,V$ be partial isometries for which $\Phi(P)=\lambda
U, \Phi(I-P)=\mu V$. We have
     \begin{equation}\label{E:krit}
     \lambda U+\mu V=I \quad \text{and} \quad
     \epsilon \lambda U+\delta \mu V \in \mathbb C \mathscr{U}(\cH)
     \quad (|\epsilon|=|\delta|=1).
     \end{equation}
Since $P\neq 0, I$, it follows that $\lambda, \mu >0$.
Choose different $\epsilon$ and $\delta$ with $|\epsilon|=|\delta|=1$.
Since by \eqref{E:krit} it follows that the operator
     \[
     \delta I+(\epsilon -\delta )\lambda U=
     \epsilon \lambda U+\delta (I-\lambda U)=
     \epsilon \lambda U+\delta \mu V
     \]
is normal, we obtain that $U$ and then that $V$ are both normal partial
isometries. Therefore, $U$ has a matrix representation
     \[
     U=\left[ \matrix U_0 & 0\\
                      0   & 0\endmatrix \right]
     \]
where $U_0$ is unitary on a proper closed linear subspace $\cH_0$ of
$\cH$. In accordance
with \eqref{E:krit}, we have the following matrix representation of $V$
     \[
     V=\left[ \matrix (I-\lambda U_0)/\mu & 0\\
                      0   & I/\mu            \endmatrix \right].
     \]
Using the characteristic property $VV^*V=V$ of partial
isometries, we get that $\mu=1$ and, by symmetry, that
$\lambda=1$. Taking the matrix
representations above into account, it is easy to see
that $I-U_0$ is a normal partial isometry and that
$\epsilon U_0+\delta (I-U_0)$ is a scalar multiple of a unitary operator
for every $\epsilon, \delta\in \mathbb C$ with $|\epsilon|=|\delta|=1$.
Since $I-U_0$ is a normal partial isometry, the spectrum of $U_0$ must
consist of such numbers
$c$ of modulus one, for which either $1-c$ has modulus one or $1-c=0$.
This gives us that $\sigma(U_0) \subset \{ 1, e^{i\pi/3},
e^{-i\pi/3}\}$.
Let $P_1,P_2,P_3$ denote the projections onto the subspaces
$\ker (U_0-I), \ker (U_0-e^{i\pi/3}I), \ker(U_0-e^{-i\pi/3}I)$
of $\cH_0$, respectively. We assert that two of
the operators $P_1, P_2, P_3$ are necessarily zero. In fact, if for
example, $P_2,P_3\neq 0$, then it follows from the second property
in \eqref{E:krit} that
     \[
     |\epsilon e^{i\pi/3}+\delta e^{-i\pi/3}|=
     |\epsilon e^{-i\pi/3}+\delta e^{i\pi/3}|
     \]
for every $\epsilon, \delta$ of modulus one. But this is an obvious
contradiction. The other cases can be treated in a similar way.
Therefore, we have $\Phi(P)=U\in \{ 1, e^{i\pi/3}, e^{-i\pi/3}\}
\mathscr{P}(\cH)$ for every projection $P$ on $\cH$.
Now, let $P$ be a projection having infinite
rank and infinite corank. Since in this case $P$ is unitarily equivalent
to $I-P$,
it follows that $P$ and $I-P$ can be connected by a continuous curve
within the set
of projections. Consequently, we obtain that $\Phi(P)$ and $\Phi(I-P)$
have the same
nonzero eigenvalue. Since $\Phi(I-P)=I-\Phi(P)$, it follows that this
eigenvalue is 1. Thus we obtain that $\Phi(P)$ is a projection.
If $P$ is a finite rank projection, then $P$ is the difference of two
projections having infinite rank and corank. Then we obtain that
$\Phi(P)$
is the difference of two projections and consequently $\Phi(P)$ is
self-adjoint. On the other hand, we have $\Phi(P)\in \{ 1, e^{i\pi/3},
e^{-i\pi/3}\} \mathscr{P}(\cH)$. These result in
$\Phi(P)\in \mathscr{P}(\cH)$ and we deduce that $\Phi$
sends every projection to a projection.
It now follows that $\Phi$ is a Jordan *-endomorphism of $\cB(\cH)$.
Since, by our condition, the range of $\Phi$ contains a rank-one
operator and an operator with dense range, using \cite[Theorem
1]{MolStud1} again, we
infer that $\Phi$ is either a *-automorphism or a *-antiautomorphism of
$\cB(\cH)$. In both cases we obtain that $\Phi$ is a surjective isometry
of $\cB(\cH)$ and this completes the proof.
\end{proof}

\begin{lemma}\label{L:unisubspace}
If $\mathscr M\subset \mathbb C \mathscr{U}(\cH)$ is a linear subspace,
then $\mathscr M$ is either 1-dimensional or 0-dimensional.
\end{lemma}

\begin{proof}
Suppose on the contrary that $\mathscr M$ is at least 2-dimensional.
Without serious loss of generality we may assume that there are linearly
independent unitaries $U,V$ in $\mathscr M$ such that
     \begin{equation}\label{E:subs1}
     I=\lambda U+\mu V
     \end{equation}
holds true for some $\lambda, \mu\in \mathbb C$. Clearly, either
$\lambda \neq 0$ or $\mu \neq 0$. Suppose that $\lambda \neq 0$. We
have
     \[
     (I-\lambda U)(I-\bar{\lambda} U^*)=|\mu |^2VV^*=|\mu|^2I
     \]
which implies
     \[
     I-2 \Re \lambda U +|\lambda|^2I=|\mu|^2I.
     \]
This means that for the unitary operator $U'=(\lambda/|\lambda|) U$ we
have
     \[
     U'+{U'}^*=2\Re U'=cI
     \]
for some real constant $c$. Multiplying this equality by $U'$, we arrive
at the equality
     \[
     {U'}^2+I=cU'.
     \]
This implies that the elements of the spectrum of $U'$ are roots of a
polynomial of degree 2. Conseqently, the spectrum of $U$ has at most two
elements. From the
original equation \eqref{E:subs1} it now follows that $U,V$ have matrix
representations
     \[
     U=\left[ \matrix aI & 0\\
                      0  & bI\endmatrix \right] , \quad
     V=\left[ \matrix cI & 0\\
                      0  & dI\endmatrix \right],
     \]
where $a,b,c,d \in \mathbb C$ are of modulus 1.
Using the condition that every  linear combination of $U$ and $V$ is a
scalar multiple of a unitary operator, we conclude that $|\epsilon
a+\delta
c|= |\epsilon b+\delta d|$ holds true for every $\epsilon, \delta \in
\mathbb C$. This readily implies that the angle between $a$ and $c$ is
the same as the angle between $b$ and $d$. Obviously, we obtain that
$U,V$ are linearly dependent which is a contradiction.
\end{proof}

\begin{lemma}\label{L:izgi}
Let $X$ be a locally compact Hausdorff space.
Let $\mathscr M\subset \mathbb C^X$ be a linear subspace containing
a nowhere vanishing function $f_0\in \mathscr M$ and having the
property that $|f|\in C_0(X)$ for every $f\in \mathscr M$.
Then there is a function $t:X \to \mathbb C$ of modulus one such that
$t\mathscr M\subset C_0(X)$.
\end{lemma}

\begin{proof}
We know that the function $|f+f_0|^2-|f|^2-|f_0|^2$ is continuous for
every $f\in \mathscr M$. This gives us that $f\overline{f_0}$ is continuous
for every $f\in \mathscr M$. Let
$t=|f_0|/f_0$. Then we have $|t|=1$ and the function
$(tf)|f_0|=(tf)(\overline{tf_0})=f\overline{f_0}$
is continuous. Consequently, we obtain $tf\in C_0(X)$.
\end{proof}

For the proof of Theorem~\ref{T:convref} we recall the well-known
Banach-Stone theorem stating that the surjective
isometries of the function algebra $C_0(X)$ are all of the form
$f \mapsto \tau \cdot f\circ \vp$, where $\tau:X \to \mathbb C$ is a
continuous function of modulus one and $\vp:X \to X$ is a homeomorphism.

\begin{proof}[Proof of Theorem~\ref{T:convref}]
Let $\Phi: \Cz\to \Cz$ be a local surjective isometry.
Pick a function $f\in C_0(X)$ and a point $x\in X$, and
consider the linear map $\Psi:\cB(\cH) \to \cB(\cH)$
     \[
     \Psi: A \longmapsto \Phi(fA)(x).
     \]
It follows from Theorem~\ref{T:formautiso} that for every $A\in
\cB(\cH)$
there exist a number $\lambda_A$ and a surjective isometry
$\tau_A\in \Iso(\cB(\cH))$ such that
$\Psi(A)=\lambda_A \tau_A(A)$. By Lemma~\ref{L:locisom} we infer that
there exist a nonnegative number $\lambda_{f,x}$ and a surjective
linear isometry $\tau_{f,x}\in \Iso(\cB(\cH))$ for which
     \begin{equation}\label{E:iso3}
     \Phi(fA)(x)=\lambda_{f,x} \tau_{f,x}(A)
     \end{equation}
holds true for every $f\in C_0(X)$, $A\in
\cB(\cH)$ and $x\in X$.
Now, let $U\in \cB(\cH)$ be a unitary operator and $x\in X$.
The linear map
     \[
     f \longmapsto \Phi(fU)(x)
     \]
maps $C_0(X)$ into $\mathbb C \mathscr U(\cH)$. Since the range of this
map is a linear subspace, by Lemma~\ref{L:unisubspace} we infer that it
is either 1-dimensional or 0-dimensional. Thus there is a linear
functional
$F_{U,x}:C_0(X)\to \mathbb C$ and a unitary operator $[\tau(x)](U)$
such that
     \begin{equation*}
     \Phi(fU)(x)=F_{U,x}(f)[\tau(x)](U) \qquad (f\in C_0(X),U\in
     \mathscr{U}(\cH), x\in X).
     \end{equation*}
Clearly, the map $F_U: C_0(X) \to \mathbb C^X$ defined by
$F_U(f)(x)=F_{U,x}(f)$ is linear and we have
     \begin{equation}\label{E:iso1}
     \Phi(fU)(x)=F_U(f)(x)[\tau(x)](U) \qquad (f\in C_0(X),U\in
     \mathscr{U}(\cH), x\in X).
     \end{equation}
Since $\Phi$ is a local surjective isometry of $\Cz$, it follows from
Theorem~\ref{T:formautiso} that for every $f\in C_0(X)$ there exist
a strongly continuous function $\tau_{f,U}:X \to \Iso(\cB(\cH))$ and a
homeomorphism $\vp_{f,U}:X\to X$ such that
     \begin{equation}\label{E:izgi}
     \Phi(fU)(x)=f(\vp_{f,U}(x))[\tau_{f,U}(x)](U) \qquad (x\in X).
     \end{equation}
Apparently, we have $|F_U(f)|=|f|\circ \vp_{f,U}$. Because of the
$\sigma$-compactness of
$X$, it is a quite easy consequence of Uryson's lemma that there exists
a strictly positive function in $C_0(X)$. Therefore, the range
of $F_U$ contains a nowhere vanishing function
and has the property that the absolute value of every function
belonging to this range is continuous. By Lemma~\ref{L:izgi}, there
exists a function $t:X \to \mathbb C$ of modulus one such that the
functions $tF_U(f)$
are all continuous $(f\in C_0(X))$. Consequently, we may suppose that
the map $F_U$ in \eqref{E:iso1} maps $C_0(X)$ into itself. Comparing
\eqref{E:iso1} and \eqref{E:izgi} we have
     \begin{equation}\label{E:izgi2}
     F_U(f)(x)[\tau(x)](U)=f(\vp_{f,U}(x))[\tau_{f,U}(x)](U) \qquad
     (x\in X).
     \end{equation}
If $f\in C_0(X)$ is a nowhere vanishing function, then by the continuity
of the functions $F_U(f),f\circ \vp_{f,U}$ and $[\tau_{f,U}(.)](U)$, it
follows
that $[\tau(.)](U)$ is also continuous. From \eqref{E:izgi2} we have
     \[
     F_U(f)=f(\vp_{f,U}(x))[\tau_{f,U}(x)](U)[\tau(x)](U)^* \qquad (x\in
     X).
     \]
In particular, this implies that the function
     \[
     x \longmapsto [\tau_{f,U}(x)](U)[\tau(x)](U)^*
     \]
can be considered as a continuous scalar valued function of modulus one.
Hence, $F_U$ is a local surjective isometry of $C_0(X)$. By our assumption
this means that $F_U$ is a surjective isometry, i.e. there exist a
continuous function $t_U:X\to \mathbb C$ of modulus one and a
homeomorphism $\vp_U:X\to X$ such that $F_U(f)=t_U \cdot f\circ \vp_U$
$(f\in C_0(X), U\in \mathscr{U}(\cH))$. Having a look at
\eqref{E:iso1}, it is obvious that we may suppose that $\Phi$ satisfies
     \begin{equation*}
     \Phi(fU)(x)=f(\vp_U(x))[\tau(x)](U) \qquad (f\in C_0(X),U\in
     \mathscr{U}(\cH), x\in X),
     \end{equation*}
where $[\tau(x)](U)$ is unitary.
If $f\in C_0(X)$ is nonnegative, we see from \eqref{E:iso3} that
     \[
     f(\vp_U(x))=\lambda_{f,x}=f(\vp_I(x))
     \]
and
     \[
     [\tau(x)](U)=\tau_{f,x}(U) \qquad (U\in \mathscr{U}(\cH), x\in X).
     \]
This verifies the existence of a homeomorphism $\vp$ of $X$ and,
due to the fact that every operator in $\cB(\cH)$ is a linear
combination of unitaries, the existence of a function $\tau:X \to
\Iso(\cB(\cH))$ for which
     \begin{equation*}
     \Phi(fU)(x)=f(\vp(x))[\tau(x)](U) \qquad (U\in \mathscr{U}(\cH),
     x\in X)
     \end{equation*}
holds true for every nonnegative function $f\in C_0(X)$. Since
every function in $C_0(X)$ is the linear combination of nonnegative
functions in $C_0(X)$, we finally obtain that
     \begin{equation*}\label{E:iso7}
     \Phi(fA)(x)=f(\vp(x))[\tau(x)](A) \qquad (f\in C_0(X), A\in
     \cB(\cH), x\in X).
     \end{equation*}
Referring to the fact once again that the linear span of the elementary
tensors $fA$ $(f\in C_0(X), A\in \cB(\cH))$ is dense in $C_0(X,
\cB(\cH))$, we arrive at the form
     \begin{equation*}
     \Phi(f)(x)=[\tau(x)](f(\vp(x))) \qquad (f\in C_0(X, \cB(\cH)),
     x\in X).
     \end{equation*}
By Theorem~\ref{T:formautiso}, the proof is complete.
\end{proof}

We now turn to the proof of Theorem~\ref{T:funcref}.
The next result describes the form of local surjective isometries of
the function algebra $C_0(X)$.

\begin{lemma}\label{L:formlociso}
Let $X$ be a first countable locally compact Hausdorff space. Let
$F:C_0(X) \to C_0(X)$ be a local surjective isometry.
Then there exist a continuous function $t:X \to \mathbb C$ of modulus
one and a homeomorphism $g$ of $X$ onto a subspace of $X$ so that
     \begin{equation}\label{E:locsurj}
     F(f)\circ g=t\cdot f \qquad (f\in C_0(X)).
     \end{equation}
\end{lemma}

\begin{proof}
By Banach-Stone theorem on the form of surjective linear isometries of
$C_0(X)$ it follows that for every $f\in C_0(X)$ there exist a
homeomorphism $\vp_f :X \to X$ and a continuous function $\tau_f:X\to
\mathbb C$ of modulus one such that
    \begin{equation}\label{E:loki}
    F(f)=\tau_f \cdot f\circ \vp_f.
    \end{equation}
For any $x\in X$ let $\mathcal S_x$ denote the set of all functions
$p\in C_0(X)$ which map into the interval $[0,1]$, $p(x)=1$ and
$p(y)<1$ for every $x\neq y\in X$. By Uryson's lemma and the first
countability of $X$, it is easy to verify that $\mathcal S_x$ is
nonempty. Let $p,p'\in \mathcal S_x$. By \eqref{E:loki} there exist
$y,y'\in X$ for which $|F(p)|\in \mathcal S_y, |F(p')|\in \mathcal
S_{y'}$. Similarly, since $(p+p')/2\in \mathcal S_x$, there is
a point $y''\in X$ for which $|F((p+p')/2)|\in \mathcal S_{y''}$.
Apparently, we have $y=y'$ and $F(p)(y)=F(p')(y')$. This shows
that there are functions $t:X \to \mathbb C$ and $g:X\to X$ such that
     \begin{equation}\label{E:maxi}
     t(x)=F(p)(g(x))
     \end{equation}
holds true for every $x\in X$ and $p\in \mathcal S_x$. Clearly,
$|t(x)|=1$. Pick $x\in X$. It is easy to see that for any
strictly positive function
$f\in C_0(X)$ with $f(x)=1$ we have a function $p\in \mathcal S_x$ such
that $p(y)< f(y)$ $(x\neq y\in X)$. Now, let $f\in C_0(X)$ be an
arbitrary
nonnegative function. Then there is a positive constant $c$ for which
the function $y\mapsto c+f(x)-f(y)$ is positive. Hence, we can choose a
function $p\in \mathcal S_x$ such that $cp(y)<cp(x)+f(x)-f(y)$ $(x\neq
y\in X)$. This means that the nonnegative function $cp+f$ takes its
maximum only at $x$. By \eqref{E:maxi} we infer
     \[
     t(x)(cp(x)+f(x))=F(cp+f)(g(x)).
     \]
Clearly, we have
     \[
     t(x)(cp(x))=F(cp)(g(x)),
     \]
too. Therefore, we obtain
     \begin{equation}\label{E:formi}
     t\cdot f=F(f)\circ g
     \end{equation}
for every nonnegative $f$ and then for every function in $C_0(X)$.
We prove that $g$ is a homeomorphism of $X$ onto
the range of $g$. To see this, first observe that for every function
$p\in \mathcal S_y$ and net $(y_\alpha)$ in $X$, the condition that
$p(y_\alpha)\to 1$ implies that $y_\alpha \to y$. Let $(x_\alpha)$ be a
net in $X$ converging to $x\in X$. Pick $p\in \mathcal S_x$. Since $F$
is a local surjective isometry, we have a homeomorphism $\vp$ of $X$ for
which
     \[
     p=|t\cdot p|=|F(p)\circ g|=p\circ \vp \circ g.
     \]
Since this implies that $p(\vp(g(x_\alpha)))\to 1$, we obtain
$\vp(g(x_\alpha)) \to x=\vp(g(x))$ and hence we have $g(x_\alpha) \to
g(x)$. So, $g$ is continuous. The injectivity of $g$ follows from
\eqref{E:formi} immediately using the fact that the nonnegative
elements of $C_0(X)$ separate the points of $X$. As for the continuity
of $g^{-1}$ and $t$, these follow
from \eqref{E:formi} again and from Uryson's lemma.
\end{proof}

Now, we are in a position to prove our last theorem.

\begin{proof}[Proof of Theorem~\ref{T:funcref}]
It is well-known that every open convex subset of $\mathbb R^n$ is
homeomorphic to the open unit ball $B$ of $\mathbb R^n$. Hence, it is
sufficient to show that the automorphism and isometry groups of $C_0(B)$
are algebraically
reflexive. Furthermore, by the form of the automorphisms and surjective
isometries of the function algebra $C_0(X)$ we are certainly done if
we prove the statement only for the isometry group. So, let $F:C_0(B)\to
C_0(B)$ be a local
surjective isometry. Then $F$ is of the form \eqref{E:locsurj}.
The only thing that we have to verify is that the function $g$ appearing
in this form is surjective.
Consider the function $f\in C_0(B)$ defined by $f(x)=1/(1+\|x\|)$.
Clearly, we may assume that $F(f)=f$. From \eqref{E:locsurj} we infer
that
     \[
     \frac{1}{1+\|x\|}=\frac{1}{1+\|g(x)\|} \qquad (x\in S).
     \]
Therefore, the continuous function $g$ maps the surface $S_r$ of the
closed ball $r\overline{B}$ $(0\leq r<1)$ into itself. It is
obvious that every proper closed subset of $S_r$ is homeomorphic
to a subset of $\mathbb R^{n-1}$.
By Borsuk-Ulam theorem we get that $g$ takes the same value at some
antipodal points of $S_r$. But this contradicts the injectivity of $g$.
Consequently, the range of $g$ contains every set $S_r$ $(0\leq r<1)$
which means that $g$ is bijective. This completes the proof.
\end{proof}

The proof of Theorem~\ref{T:funcref} shows how difficult it might be to
treat our reflexivity problem for the suspension of arbitrary
$C^*$-algebras. We mean the role of the use of Borsuk-Ulam
theorem in the above argument. To reinforce this opinion, let us
consider only
the particular case of commutative $C^*$-algebras. Let $X$ be a locally
compact Hausdorff space and suppose that the automorphism and isometry
groups of $C_0(X)$ are algebraically reflexive. If $F:C_0(\mathbb
R \times X)\to C_0(\mathbb R \times X)$ is a local surjective
isometry, then Lemma~\ref{L:formlociso} gives the form of $F$. The
problem is to verify that the function $g$ appearing in
\eqref{E:locsurj} is surjective. This would be easy if there were an
injective nonnegative function in $C_0(\mathbb R \times X)$.
Unfortunately, this is not the case even when $X$ is a singleton.
Anyway, if $n\geq 3$, there
is no injective function in $C_0(\mathbb R^n)$ at all. Therefore,
to attack the
problem of the surjectivity of $g$, we had to invent a different
approach which was the use of Borsuk-Ulam theorem. To mention another
point, it is easy to see that in general the automorphisms as well as
surjective
isometries of the tensor product $C_0(X_1)\otimes C_0(X_2)\cong
C_0(X_1\times X_2)$ have nothing to do with the automorphisms and
surjective isometries of $C_0(X_1)$ and $C_0(X_2)$, respectively.
However, according to Theorem~\ref{T:formautiso}, in the case
of the tensor product $C_0(X)\otimes \cB(\cH)$ every automorphism as
well as surjective isometry is an easily identifiable mixture of a
"functional algebraic" and an "operator algebraic" part. This
observation was of fundamental importance when verifying the result in
Corollary~\ref{C:suspref}.
These might justify the suspicion why we feel our reflexivity
problem really difficult even for the supsension of general commutative
$C^*$-algebras.



\begin{thebibliography}{99999}

\bibitem[BaMo]{BaMo}
C.J.K. Batty and L. Moln\'ar,
On topological reflexivity of the groups of *-automorphisms
and surjective isometries of $B(H)$,
\emph{Arch. Math.} \textbf{67} (1996), 415--421.

\bibitem[Bre]{Bre}
M. Bre\v sar,
Characterizations of derivations on some normed algebras with
invloution,
\emph{J. Algebra} \textbf{152} (1992), 454--462.

\bibitem[BrSe1]{BrSe1}
M. Bre\v sar and P. \v Semrl,
Mappings which preserve idempotents, local automorphisms, and
local derivations,
\emph{Canad. J. Math.} \textbf{45} (1993), 483--496.

\bibitem[BrSe2]{BrSe2}
M. Bre\v sar and P. \v Semrl,
On local automorphisms and mappings that preserve idempotents,
\emph{Studia Math.} \textbf{113} (1995), 101--108.

\bibitem[Che]{Che}
P.R. Chernoff,
Representations, automorphisms and derivations of some operator
algebras,
\emph{J. Funct. Anal.} \textbf{12} (1973), 275--289.

\bibitem[CiYo]{CiYo}
P. Civin and B. Yood,
Lie and Jordan structures in Banach algebras,
\emph{Pacific J. Math.} \textbf{15} (1965), 775--797.

\bibitem[Cri]{Cri}
R.L. Crist,
Local derivations on operator algebras,
\emph{J. Funct. Anal.} \textbf{135} (1996), 76--92.

\bibitem[DFR]{DFR}
T. Dang, Y. Friedman and B. Russo,
Affine geometric proofs of the Banach Stone theorems of Kadison
and Kaup,
\emph{Rocky Mountain J. Math.} \textbf{20} (1990), 409--428.

\bibitem[Had]{Had}
D. Hadwin,
A general view of reflexivity,
\emph{Trans. Amer. Math. Soc.} \textbf{344} (1994), 325--360.

\bibitem[Kad]{Kad}
R.V. Kadison,
Local derivations,
\emph{J. Algebra} \textbf{130} (1990), 494--509.

\bibitem[Lar]{Lar}
D.R. Larson,
Reflexivity, algebraic reflexivity and linear interpolation,
\emph{Amer. J. Math.} \textbf{110} (1988), 283--299.

\bibitem[LaSo]{LaSo}
D.R. Larson and A.R. Sourour,
Local derivations and local automorphisms of $B(X)$,
\emph{in} Proc. Sympos. Pure Math. 51, Part 2, Providence, Rhode
Island 1990, 187--194.

\bibitem[Mol1]{MolStud1}
L. Moln\'ar,
The set of automorphisms of $B(H)$ is topologically reflexive in
$B(B(H))$,
\emph{Studia Math.} \textbf{122} (1997), 183--193.

\bibitem[Mol2]{MolLond}
L. Moln\'ar,
Reflexivity of the automorphism and isometry groups of
$C^*$-algebras in BDF theory,
submitted for publication.

\bibitem[Mur]{Mur}
G.J. Murphy,
"$C^*$-algebras and Operator Theory,"
Academic Press, 1990.

\bibitem[Nai]{Nai}
M.A. Naimark,
"Normed Algebras,"
Wolters-Noordhoff Publishing, 1972.

\bibitem[Shu]{Shu}
V.S. Shul'man,
Operators preserving ideals in $C^*$-algebras,
\emph{Studia Math.} \textbf{109} (1994), 67--72.

\bibitem[ZhXi]{ZhXi}
J. Zhu and C. Xiong,
Bilocal derivations of standard operator algebras,
\emph{Proc. Amer. Math. Soc.} \textbf{125} (1997), 1367--1370.

\end{thebibliography}
\end{document}